\documentclass{article}
\usepackage{graphicx}
\usepackage{lipsum}  
\usepackage{amsmath} 
\usepackage{enumitem} 
\usepackage{amsfonts} 
\usepackage{geometry} 
\usepackage{graphicx} 
\usepackage{caption} 
\usepackage{subcaption} 
\usepackage[ruled,vlined]{algorithm2e} 
\usepackage{listings} 
\usepackage{booktabs}

\begin{document}

\title{Tumoral Angiogenesis Optimizer: A new bio-inspired based metaheuristic}
\author{
  {\small Hernández Rodríguez, Matías Ezequiel}\thanks{E-mail: \texttt{mhernandez@inidep.com}} \\ \\
  {\small Department of Mathematical Biology} \\
  {\small Instituto Nacional de Investigación y Desarrollo Pesquero} \\
  {\small Mar del Plata, Argentina} 
}

\date{} 

\maketitle

\begin{abstract}
In this article, we propose a new metaheuristic inspired by the morphogenetic cellular movements of endothelial cells (ECs) that occur during the tumor angiogenesis process. This algorithm starts with a random initial population. In each iteration, the best candidate selected as the tumor, while the other individuals in the population are treated as ECs migrating toward the tumor's direction following a coordinated dynamics through a spatial relationship between tip and follower ECs. 
This algorithm has an advantage compared to other similar optimization metaheuristics:
the model parameters are already configured according to the tumor angiogenesis phenomenon modeling, preventing researchers from initializing them with arbitrary values.
Subsequently, the algorithm is compared against well-known benchmark functions, and the results are validated through a comparative study with Particle Swarm Optimization (PSO). The results demonstrate that the algorithm is capable of providing highly competitive outcomes.
Furthermore, the proposed algorithm is applied to real-world problems (cantilever beam design, pressure vessel design, tension/compression spring and sustainable explotation renewable resource). The results showed that the proposed algorithm worked effectively in solving constrained optimization problems. The results obtained were compared with several known algorithms.
\end{abstract}


\section{Introduction}

Optimization is a broad concept that permeates various domains, from engineering design to business planning, from internet routing to environmental sustainability. Businesses seek to maximize profits while minimizing costs, engineers strive to optimize the performance of their designs while minimizing expenses, and sustainability studies aim to minimize environmental harm in resource exploitation. Nearly everything we do is, in some way, related to optimizing something. Consider, for instance, vacation planning, where we seek to maximize enjoyment while minimizing costs \cite{Diw, Gil, Gol, Hul, Len, Loz, Lupt, Set, Yan1}.

Mathematical optimization addresses these problems using mathematical tools. However, it has the drawback that many algorithms, especially gradient-based search methods, are local search techniques. Typically, these searches begin with an assumption and attempt to improve the quality of solutions. For unimodal functions, convexity ensures that the final optimal solution is also a global optimum. For multimodal objectives, the search may become trapped in a local optimum. Another limitation lies in solving optimization problems with high-dimensional search spaces; classical optimization algorithms fail to provide suitable solutions because the search space grows exponentially with problem size, rendering exact techniques impractical. Moreover, the complexities of real-world problems often prevent verifying the uniqueness, existence, and convergence conditions that mathematical methods require \cite{Kel, Yan1}.

Metaheuristics, on the other hand, are becoming powerful methods to solve many challenging optimization problems. They are known for their ability to find global optima. Classic examples include genetic algorithms, ant colony optimization (ACO), and PSO, among others. Metaheuristics find applications across science, technology, and engineering fields \cite{Hau, Mir1, Mir2, Yan1}.

In recent years, there has been a growing interest in nature-inspired, human behavior-inspired, and physics-inspired metaheuristics, such as the Moth Search Algorithm, Grey Wolf Optimizer (GWO), Gold Rush Optimizer (GRO), Bat-Inspired Algorithm (BA), Ebola optimization search algorithm (EOSA), and the Gravitational Search Algorithm (GSA) \cite{Mir0, Oye, Ras, Yang, Zol}. These algorithms address different optimization problems, but there is no one-size-fits-all algorithm that provides the best solution for all optimization problems. Some algorithms perform better than others for specific problems. Hence, the search for new heuristic optimization algorithms remains an open problem.

In this paper, we introduce a novel optimization algorithm inspired by the morphogenetic cell movements of ECs that occur during tumor angiogenesis, namely, the Tumoral Angiogenesis Optimizer (TAO). This article is structured as follows: In Section 2, we describe an agent-based model that explains the behavior of endothelial cells during angiogenesis. This model, along with the PSO algorithm, inspired our optimization algorithm, which is detailed in Section 3. In Section 4, a comparative study with the PSO algorithm is presented, considering test functions. In the same section, constrained optimization problems are addressed (Rosenbrock function constrained with a cubic and a line, cantilever beam design, pressure vessel design, tension/compression spring and sustainable explotation renewable resource). Finally, conclusions and possible new researches are presented in Section 5.

\section{Angiogenic cell movements mathematical model}

Morphogenetic cell movements generate diverse tissue and organ shapes, and questions arise regarding whether these movements share common principles and how cells coordinate behaviors. Angiogenesis, a morphogenetic cell movement, involves the emergence of new vascular networks. Vascular ECs work collectively to form dendrite structures, with several molecular players identified. However, the cellular mechanisms underlying angiogenesis remain largely unknown. Understanding these processes would bridge the gap between molecules and angiogenic morphogenesis.

To gain deeper insights into morphogenesis, a group of researchers developed a system that combines time-lapse imaging with computer-assisted quantitative analysis. This approach allowed for a thorough investigation of the behaviors of ECs driving angiogenic morphogenesis in an in vitro model \cite{Ari}. The discoveries revealed that EC behaviors are considerably more dynamic and intricate than previously believed, with individual ECs frequently changing their positions, including instances of tip cell overtaking \cite{Ari}. The phenomenon of dynamic tip cell overtaking had also been reported by another research group \cite{Jak}.
 These revelations led the researchers to ponder the following question: How are the movements of individual ECs integrated into the highly dynamic and complex multicellular process that culminates in the formation of ordered architectures?
Mathematical and computational modeling strategies have proven invaluable for shedding light on the biological intricacies underlying angiogenic morphogenesis, especially when employed alongside quantitative experimental approaches \cite{Cha, Pei, Sci, Tys}. Over time, a variety of models, encompassing continuum, discrete, and hybrid approaches, have been developed to explore different facets of angiogenesis across various biological scales \cite{Arak, Bai, Pet}. Recent advancements have introduced cell-based models, such as cellular potts and agent-based models, designed to uncover the biological implications of angiogenesis predictively. These models enable the dissection of the molecular and cellular mechanisms involved in processes like sprouting \cite{And, Bau, Merks, Szabo} and cell rearrangement \cite{Bent}.

In \cite{Sug}, one-dimensional an agent-based model was developed in order to simulate the behaviors of ECs during the elongation of blood vessels. In this models, individual ECs are represented as agents aligned along the axis of vessel elongation, which forms an emerging sprout in the vascular network. 
Each cell (agent) behaves autonomously in accordance with a set of specific rules:
For each step, $t,$ each agent, $i$, has a position, $x_i(t)$, a cell migration speed, $v_i(t)$ ($v_i = v_1$ or $v_2$, where $v_2 < v_1$), and a cell migration direction, $D_i(t)$ ($D = +1$ (anterograde) or $D_i = -1$ (retrograde)). For each step, $v_i$ and $D_i$ satisfy the following rules:

\begin{enumerate}
    \item Speed transition rule: If $v_i(t) = v_1$, it can change to $v_i(t+1) = v_2$ with a probability $p$. If $v_i(t) = v_2$, it can change to $v_i(t+1) = v_1$ with a probability $s$. In the absence of these conditions, motility remains unchanged, i.e., $v_i(t+1) = v_i(t)$.
    \item Direction transition rule: If $D_i(t) = -1$, then with a probability $r$, motility changes to $D_i(t+1) = +1$. If at the last step $D_i(t) = +1$, then with a probability $s$, motility changes to $D_i(t+1) = -1$. In the absence of these conditions, the direction remains unchanged, i.e., $D_i(t+1) = D_i(t)$.
\end{enumerate}

Furthermore, from the set $X_t = \{x_j(t)\}_j$, we consider $\max X_t$ and $\max\left\{\min\{x_j,x_k\}\right\}_{j\neq k}$, which are respectively the largest and second largest elements among the set $X_t$. We define $(i, t)$ satisfying a tip EC restriction condition if $\max X_t-\max\left\{\min\{x_j,x_k\}\right\}_{j\neq k}>d,$ where $d$ is the tip-follower threshold for restriction of tip EC movement. If the tip EC restriction condition is satisfied, $v_i(t) = v_2$ is adopted as the tip EC speed. This models the fact that more mature ECs play a crucial role in controlling or regulating the mobility of tip ECs.

The agents update their position based on the following equation:

\begin{equation} \label{eq1}
x_i(t+1) = x_i(t) + v_i(t+1)D_i(t+1),
\end{equation}
together with the tip EC restriction condition.

\section{Tumoral Angiogenesis Optimizer}
In this section, we explain our optimization algorithm. As we have just seen, during the angiogenesis process, there are local rules that ECs follow, but there is also global communication through which more mature ECs can regulate the migration speed of tip ECs during the formation of new blood vessels. These cells migrate towards the direction where there is a gradient of growth factors, such as vascular endothelial growth factor (VEGF), which attracts cells to the tissue where vascularization is needed. In the context of tumor angiogenesis, tumors can secrete a signaling protein to stimulate the formation of new blood vessels that grow towards the tumor to supply it with oxygen and nutrients, which is essential for its growth and survival.

This global communication among ECs reminds us of the PSO algorithm, which was initially introduced by James Kennedy and Russell Eberhart in 1995 \cite{Ken}, where global communication plays a central role and inspires us to formulate our algorithm as follows:

TAO is initialized with a population of random solutions within the problem's search space. The initial migration speeds are all set to $v_1$, and the initial migration directions are all set to 1, i.e., $v_i(0) = v_1$ and $D_i(0) = 1$ for each individual $i$ in the initial population.

In each iteration, the algorithm seeks optima by updating the individuals in the population, and the solution with the best fitness is referred to as ``the tumor", and the other solutions, referred to as ``cells," migrate through the problem's search space towards the tumor following the following dynamics:

\begin{equation} \label{eq2}
x_i(t+1) = x_i(t) + v_i(t+1)D_i(t+1)(tumor(t) - x_i(t)) + \gamma^t r,
\end{equation}
where $x_i(t)$ and $tumor(t)$ are, respectively, the position of cell $i$ and the best solution in iteration $t$. The migration speed, $v_i(t),$ and migration direction, $D_i(t),$ follow the rules 1 and 2, explained earlier, with respect to the direction from cell $x_i(t)$ to $tumor(t)$. Additionally, for each cell, the distance it has traveled throughout the algorithm is recorded, and this information is used to check if the restriction of tip EC movement is satisfied. When the tumor is renewed, the distances of the other cells are reset to 0, meaning that the capillaries are pruned. To simulate the branching of blood vessels, we used a random vector $r$, which allows us to modify the direction of $tumor(t) - x_i(t)$ above its normal plane. This provides diversity among the cells and ensures good exploration. Parameter $\gamma$ is a learning parameter, which is chosen in $[0.5,1).$ In this paper we used $\gamma = 0.7$

\begin{algorithm}
\caption{Tumoral Angiogenesis Optimizer Algorithm}
\label{alg_tao}  
\KwData{Objective function, $f(x);$ model parameters, $v_1$, $v_2$, $p$, $s$, $q,$ $d;$ population size, $N_{pop};$ max number of iterations, $N_{max};$ and search region, $\Omega$}

\KwResult{Approximate optimum, $x^+,$ and $f(x^+)$}

Initialize cell population randomly inside the search region\;
Initialize cellular migration speeds to $v_1$\;
Initialize cellular migration directions to 1\;
Initialize length associated with each cell to 0\;
Calculate the fitness of each search cell\;
$Tumor(0)\Leftarrow $ the best search cell\;

\While{$t<N_{max}$}{
    Check if tip EC restriction condition is satisfied\;
    \ForEach{cell $\neq$ Tumor(t)}{
        Check rules 1. and 2.\;
        Update the positions of the cells from equation (\ref{eq2})\;
        \If{$f(cell)<f(Tumor(t))$}{
            		$Tumor(t) \Leftarrow cell$\; 
            		Set the lengths associated with the cells equal to 0\; 
            		}     
    }
\Return{$x^+=Tumor$ and $f(x^+).$}
}

\end{algorithm}

Parameters $v_1$, $v_2$, $p$, $s$, $q$, and $d$ already set according \cite{Sug}, which provides an advantage compared to other metaheuristics where calibrating the involved parameters can be a very challenging task. From the pseudocode \ref{alg_tao}, it is relatively straightforward to implement TAO algorithm in any programming language. Table \ref{Tab0} shows values used in this paper.

\begin{table}[ht]
	\noindent\begin{tabular*}{\textwidth}
	{c @{\extracolsep{\fill}}ccccccc }
	\toprule
	$v_1$& $v_2$& $p$& $q$& $r$& $s$& $d$\\
	\hline
	5.332 & 0.938& 0.0416891& 0.234& 0.194& 0.240& 55\\
	\hline

	\end{tabular*}
\caption{Parameter values used in this paper.}
\label{Tab0}   
\end{table}

\section{Validation and comparation}

\subsection{Unconstrained optimization problems}
In order to evaluate the performance of our algorithm, we applied it to 6 standard benchmark functions. The function descriptions are detailed in Table \ref{Tab1}.

Given the stochastic nature of the meta-heuristic algorithms, their performance cannot be accurately assessed based on a single run. In order to evaluate the approach comprehensively, multiple trials with independently initialized populations are conducted. Consequently, TAO algorithm was run 50 times on each benchmark function, with a population size of 100 cells and maximum number of 500 iterations employed for both low- and high-dimensional problems.

\begin{table}[ht]
	\noindent\begin{tabular*}{\textwidth}
	{l @{\extracolsep{\fill}}ccccc }
	\toprule
	Function& Name& Dimension& Range& Minimum\\
	\hline
 	\footnotesize$F_1=\displaystyle\sum_{i=1}^{n}x_i^2$& Sphere& 20& $[-100,100]^n$& 0\\
 	\footnotesize$F_2=\displaystyle\sum_{i=1}^{n}100\left(x_i-x_{i-1}^2\right)+(x_{i-1}-1)^2$& Rosenbrock& 10& $[-30,30]^n$& 0\\
 	\footnotesize$F_3=x^2+y^2+25\left(\sin^2(x)+\sin^2(y)\right)$& Eggcrate& 2& $[-2\pi,2\pi]^2$& 0\\
 	\footnotesize$F_4=\displaystyle\sum_{i=1}^{n}(x_i+0.5)^2$& Step& 30& $[-5.12,5.12]^n$& 0\\
 	\footnotesize$F_5=10n+\displaystyle\sum_{i=1}^{n}\left[x_i^2-10\cos(2\pi x_i\right]$& Rastrigin& 10& $[-5.12,5.12]^n$& 0\\
 	\footnotesize$F_6=-\displaystyle\sum_{i=1}^{n}\sin(x_i)\left[\sin\left(\frac{ix_i^2}{\pi}\right)\right]^{20}$& Michalewicz& 5& $[0,\pi]^n$& -4.6877...\\
\footnotesize$F_7=\displaystyle\sum_{i=1}^{n}ix_i^2$& Sum Squares& 30& $[-10,10]^n$& 0\\
	\hline
	\end{tabular*}
\caption{Classical benchmark problems.}
\label{Tab1}   
\end{table}

In order to compare the performance of TAO algorithm we used standard PSO. Let's remember that the PSO algorithm has the following equations to update the velocities and positions of the particles for each particle $i$:

\begin{equation} \label{eq3}
\begin{cases}
	v_i(t+1) = wv_i(t) + c_1r_1(p_{i_{best}}-x_i(t)) + c_2r_2(g_{best}-x_i(t))\\
	x_i(t+1) = x_i(t) + v_i(t+1),
\end{cases}
\end{equation}
where $v_i(t),$ $x_i(t)$ and $p_{i_{best}}$ are, respectively, the velocity, position and personal best value of particle $i$ at time $t.$ Furthermore, $w$ is the inertial weight, $c_1$ shows the individual coefficient, $c_2$ signifies the social coefficient, $r_1,$ $r_2$ are random numbers uniformly distributed in $[0, 1],$ $g_{best}$ is the global best, and shows the best solution found by all particles (entire swarm) until $t$th iteration. Usually $c_1 = c_2 \in [1.8, 2.0].$ We used $c_1=c_2=2$ \cite{Shi}.

The inertia weight $w$ is inspired by the PSO algorithm and is calculated using the equation:

\begin{equation} \label{eq4}
w = w_{max} - \frac{w_{max} - w_{min}}{t_{max}}t,
\end{equation}
where $t$ is the current iteration step, $t_{max}$ is the maximum iteration step, $w_{max}$ is the maximum inertia weight, and $w_{min}$ is the minimum inertia weight. Usually, $w_{max} = 0.9$ and $w_{min} = 0.4$ \cite{Shi, Zha}.

The initial velocities in the PSO algorithm are set to zero, as empirical evidence has shown this to be the most effective approach \cite{Enge}.

\begin{figure}[ht]
  \begin{subfigure}[b]{0.4\textwidth}
    \includegraphics[width=\textwidth, height=\textwidth]{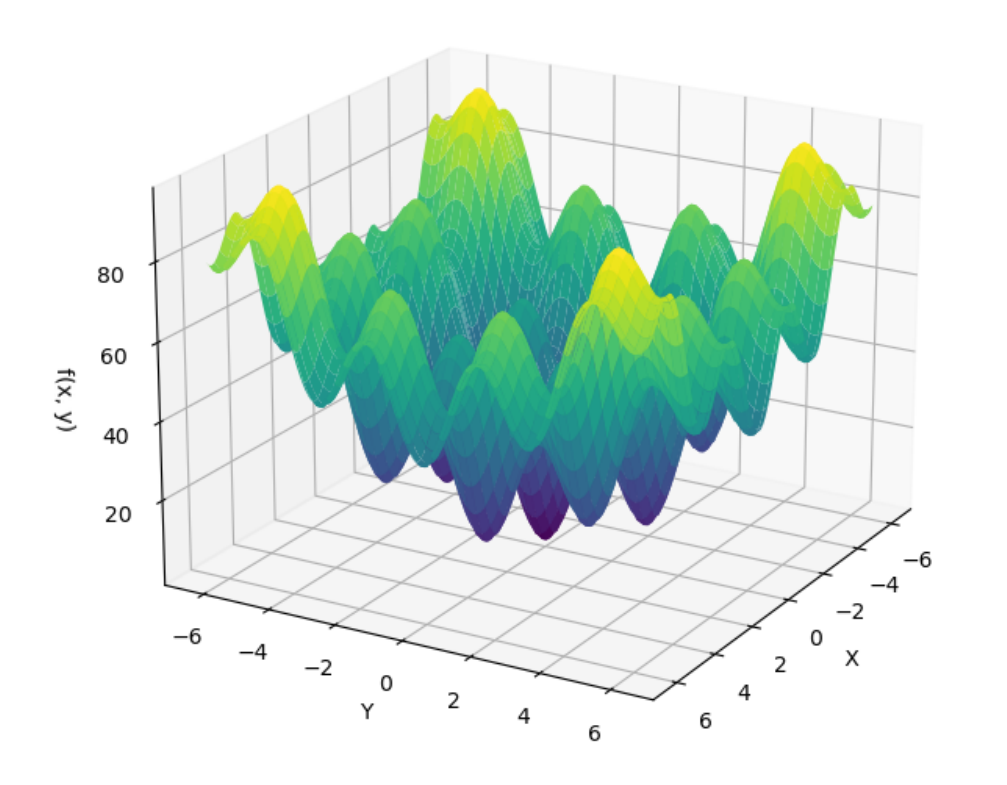}
    \caption{Eggcrate function.}
    \label{fig:f1}
  \end{subfigure}
  \hfill
  \begin{subfigure}[b]{0.4\textwidth}
    \includegraphics[width=\textwidth, height=\textwidth]{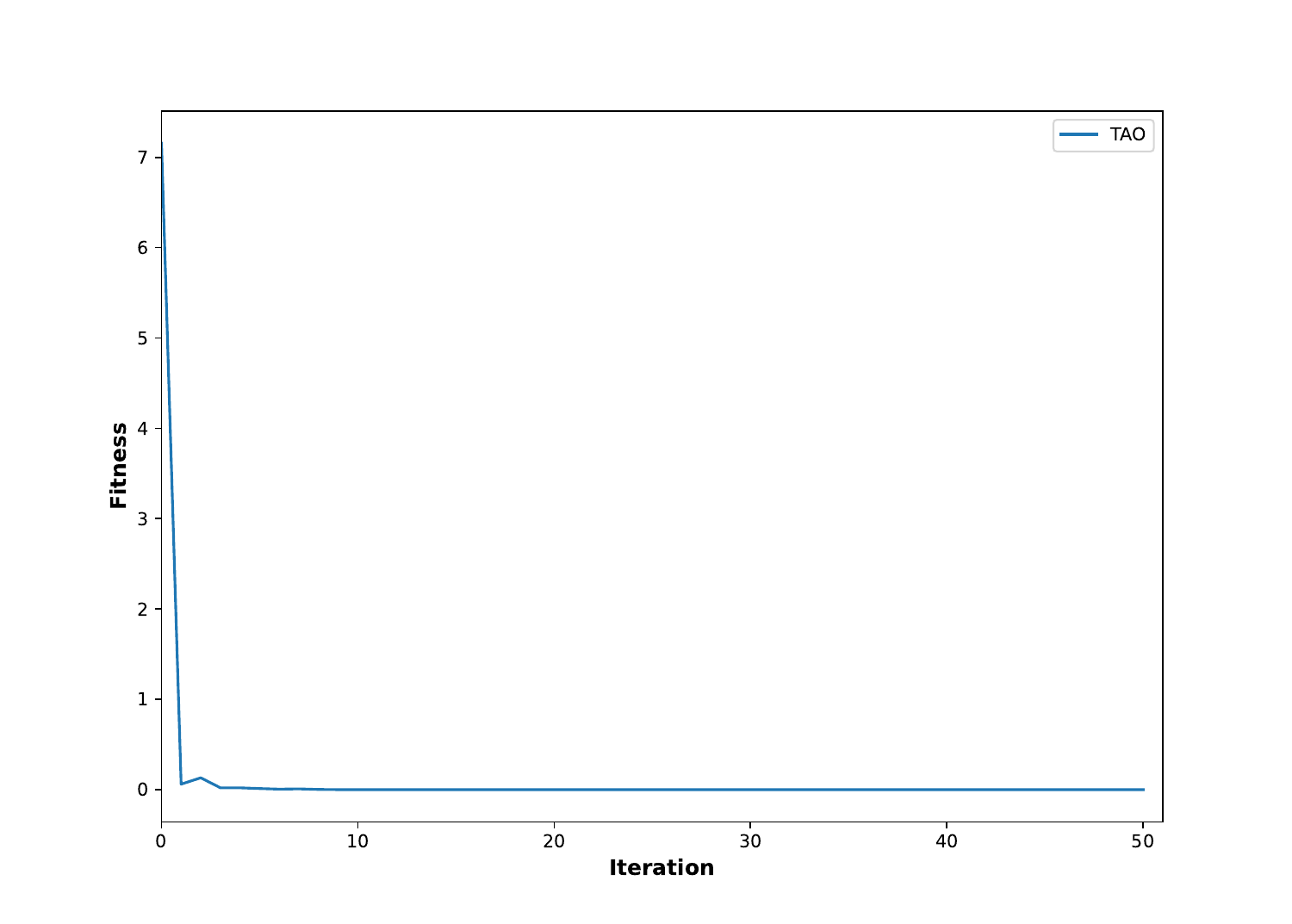}
    \caption{Convergence curve of TAO algorithm.}
    \label{fig:f2}
  \end{subfigure}
  \caption{A single run of the algorithm yielded the following approximate solution: $x^+$ = (1.3119482e-9, 3.9670538e-9), with $f(x^+)$ = 4.5444231e-16.}
  \label{figura1}
\end{figure}

\begin{table}[ht]
	\noindent\begin{tabular*}{\textwidth}
	{c @{\extracolsep{\fill}}cccc }
	\toprule
	Function& Parameters& PSO& TAO\\
	\hline
		 & Beast& 0.0416891& 1e-07\\
 	$F_1$& Mean& 1.7180551& 1.0434957\\
 		 & Standard deviation& 1.4389095& 0.8459855\\
	\hline
		 & Beast& 2.7253815& 0.0058221\\
 	$F_2$& Mean& 8.5999977& 6.9607255\\
 		 & Standard deviation& 4.6557892& 5.0302997\\
	\hline
		 & Beast& 0.0000000& 0.0000000\\
 	$F_3$& Mean& 0.0000000& 0.0000000\\
 		 & Standard deviation& 0.0000000& 0.0000000\\
	\hline
		 & Beast& 0.306966& 5.5e-06\\
 	$F_4$& Mean& 0.7906888& 0.0010151\\
 		 & Standard deviation& 0.2765826& 0.0015117\\
	\hline
		 & Beast& 1.6399746& 0.9899181\\
 	$F_5$& Mean& 7.3307733& 8.4788214\\
 		 & Standard deviation& 3.5371322& 5.4759516\\
	\hline
		 & Beast& -4.8693888& -4.6458954\\
 	$F_6$& Mean& 4.6169278& -3.9887314\\
 		 & Standard deviation& 0.141926& 0.5212977\\
	\hline
		 & Beast& 0.6684421& 0.181605\\
 	$F_7$& Mean& 8.0236676& 1.9926416\\
 		 & Standard deviation& 6.2222502& 1.7023876\\
	\hline
	\end{tabular*}
\caption{Statistical results of benchmark functions.}
\label{Tab2}   
\end{table}

Statistical results of the performance of TAO and PSO algorithms on the benchmark problems presented in Table \ref{Tab2} indicate that the first algorithm outperforms the second in all cases.

\subsection{Constrained optimization problems}
Constrained optimization problems, encompassing both equality and inequality constraints, play a pivotal role in scientific and engineering disciplines, enabling the modeling of systems subject to physical laws, budget constraints, and operational requirements. Their applications span from structural design and resource allocation to parameter estimation and process optimization.

Optimization problems with constraints pose substantial challenges for metaheuristic algorithms. The complex nature of constraint handling necessitates innovative approaches.

Mathematically, constrained minimization problem are typically expressed as:

\begin{align} \label{eqOpt}
    & \min_{x\in\mathbb{R}^n} \quad  f(x)   \\
    & \text{Subject to:}  \nonumber \\
    & g_i(x)\leq 0, i =1,\dots, N \nonumber \\
    & h_j(x)=0, j=1,\dots, P \nonumber \\ 
    & l_h \leq x_h \leq u_h, h =1,\dots, n, \nonumber 
\end{align}
where $g_j$ are inequality constraints, $h_j$ are equality constraints, $l_i$ and $u_i$ are lower and upper bounds of $x_i,$ and $f(x)$ is the objective function that needs to be optimized subject to the constraints.

There exist five prominent constraint handling methodologies, which encompass: 1) penalty functions, 2) specialized representations and operators, 3) repair algorithms, 4) the separation of objectives and constraints, and 5) hybrid methods. Among these, penalty functions represent the most straightforward approach \cite{Coello}.

In this article, we employ TAO algorithm using the static penalty approach. Essentially, depending on the nature of the problem, as demonstrated in the application examples, we pursue two strategies. 
One strategy, the more natural approach, entails transforming problem (\ref{eqOpt}) in the following problem:

\begin{align} \label{eqOpt1}
    & \min_{x\in\mathbb{R}^n} \quad  F(x)   \\
    & l_h \leq x_h \leq u_h, h =1,\dots, n \nonumber 
\end{align}
where $F(x)=\displaystyle\sum_{i=1}^Nr_i\max\{g_i(x),0\}+\displaystyle\sum_{j=1}^P c_j|h_j(x)|,$  functions $\max\{g_i(x),0\}$ and $|h_j(x)|$ measure the extent to which the constraints are violated, and $r_i$ and $c_j$ are known as penalty parameters.

In the other strategy, 

\begin{equation}\label{eqOpt2}
F(x) =
\begin{cases}
f(x) & \text{if the solution is feasible;} \\
K - \displaystyle\sum_{i=s} \frac{K}{m} & \text{otherwise,}
\end{cases}
\end{equation}
where $s$ is the number of constraints satisfied, $m$ is the total number of (equality and inequality) constraints and $K$ is a large constant ($K=1\times 10^9$) \cite{Mora}. 

\subsubsection{Rosenbrock function constrained with a cubic and a line}

Let us first consider the problem of minimizing the Rosenbrock function subject to two inequality constraints: one is a cubic constraint, and the other is linear. This is a nonlinear constrained optimization problem with a global optimum at $x^*=(1,1),$ with $f(x^*)=0.$ 

\begin{align} \label{eq5}
    & \min_{(x,y) \in [-100,100]^2} \quad  f(x,y) = (1-x)^2 + 100(y-x^2)^2  \\
    & \text{Subject to:}  \nonumber \\
    & g_1(x):(x-1)^3-y+1 \leq 0 \nonumber \\
    & g_2(x):x+y-2 \leq 0. \nonumber 
\end{align}

In order to solve problem (\ref{eq5}), we rewrite it as follows:

\begin{align} \label{eq6}
    & \min_{(x,y) \in [-100,100]^2} \quad  F(x), 
\end{align}
where $$F(x,y) = f(x,y)+\max\{g_1(x),0\}+\max\{g_2(x),0\}.$$

TAO metaheuristic converges to the approximate solution $x^+=(0.99999941, 1.00000075),$  with an absolute error  $\epsilon\approx 10^{-7}.$ 

\subsubsection{Cantilever Beam Design Problem} 

\begin{figure}[ht]
    \centering
    \includegraphics[width=0.8\textwidth]{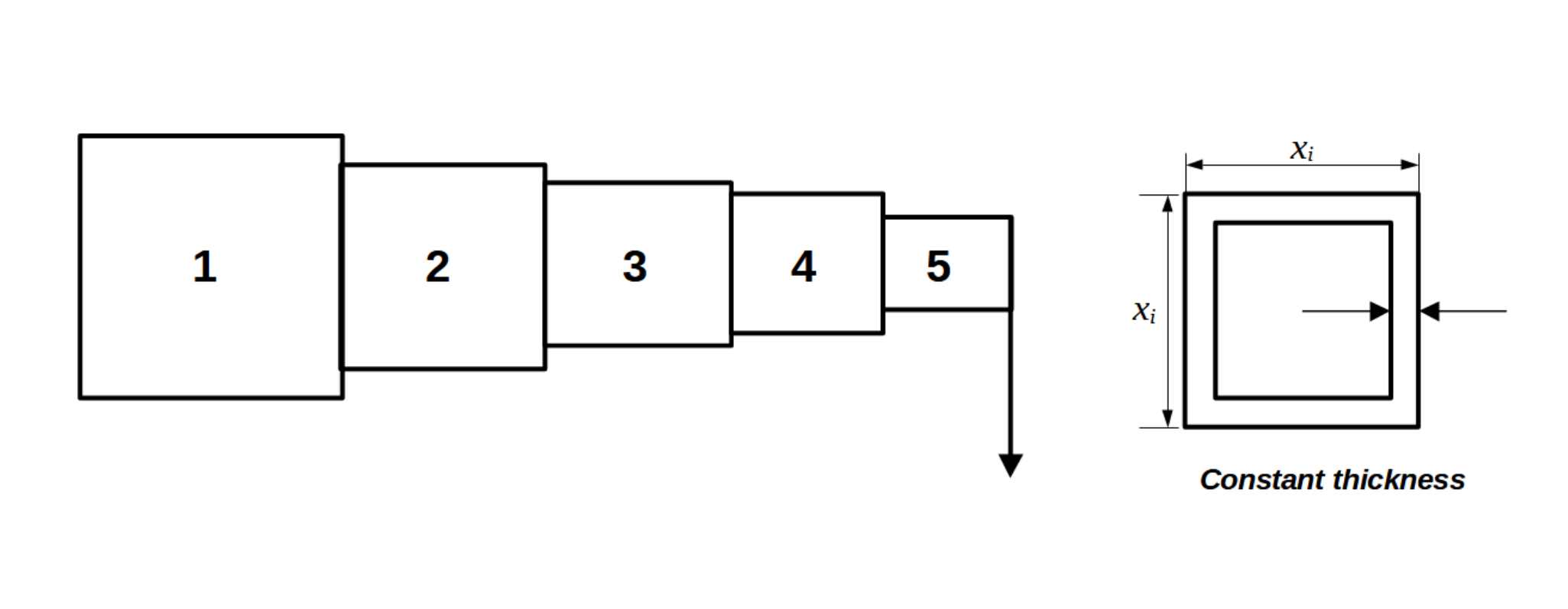}
    \caption{Cantilever Beam Design Problem.}
    \label{figura2}
\end{figure}

The Cantilever Beam Design Problem is a significant challenge in the fields of mechanics and civil engineering, primarily focused on minimizing the weight of a cantilever beam. In this problem, the beam consists of five hollow elements, each with a square cross-section. The goal is to determine the optimal dimensions of these elements while adhering to certain constraints (see figure \ref{figura2}).

The mathematical expression governing this problem and its associated constraints are represented by equation (\ref{eq7}).

\begin{align} \label{eq7}
    & \min_{0.01\leq x_1,x_2,x_3,x_4,x_5\leq 100} \quad  f(x_1,x_2,x_3,x_4,x_5) = 0.06224(x_1+x_2+x_3+x_4+x_5)  \\
    & \text{Subject to:}  \nonumber \\ 
    & g(x):\frac{61}{x_1^3} + \frac{37}{x_2^3} + \frac{19}{x_3^3} + \frac{7}{x_4^3} + \frac{1}{x_5^3} \leq 1. \nonumber 
\end{align}

In order to address problem (\ref{eq7}), we reformulate it as an unconstrained problem as follows:

\begin{align} \label{eq8}
    & \min_{0.01\leq x_1,x_2,x_3,x_4,x_5\leq 100} \quad  F(x_1,x_2,x_3,x_4,x_5), 
\end{align}
where $F(x_1,x_2,x_3,x_4,x_5)=f(x_1,x_2,x_3,x_4,x_5)+\max\{g(x)-1, 0\}.$ 
Subsequently, we apply TAO algorithm using a population of 100 individuals and a maximum of 300 iterations. Then, we obtain the approximate solution: $x^+=(6.01601588, 5.30917383, 4.49432957, 3.50147495, 2.15266534),$ with $f(x^+)=1.33652057.$

Table \ref{Tab3} lists the best solutions obtained by TAO and various methods: artificial ecosystem-based optimization (AEO), ant lion optimizer (ALO), coot optimization algorithm (COOT), cuckoo search algorithm (CS), gray prediction evolution algorithm based on accelerated even (GPEAae), hunter-prey optimizer (HPO), interactive autodidactic school (IAS), multi-verse optimizer (MVO) and symbiotic organisms search (SOS) \cite{Naru}.

The comparative outcomes are presented within  Table \ref{Tab3}. These results clearly demonstrate that the TAO algorithm we propose presents a superior solution for addressing this problem.

\begin{table}[ht]
{\small
    \noindent\begin{tabular*}{\textwidth}
    {l @{\extracolsep{\fill}}llllll }
    \toprule
    \multicolumn{1}{l}{\small Algorithm}& \multicolumn{5}{l}{\small Optimum variables}& \multicolumn{1}{c}{\small Optimum weight}\\
    \cline{2-6} 
        &$x_1$& $x_2$& $x_3$& $x_4$& $x_5$& \\
     TAO& 6.01601588& 5.30917383& 4.49432957& 3.50147495& 2.15266534& 1.33652057 \\
	 AEO& 6.02885000& 5.31652100& 4.46264900& 3.50845500& 2.15776100& 1.33996500 \\   
     ALO& 6.01812000& 5.31142000& 4.48836000& 3.49751000& 2.15832900& 1.33995000 \\
    COOT& 6.02743657& 5.33857480& 4.49048670& 3.48343700& 2.13459100& 1.33657450 \\
       CS& 6.0089000& 5.30490000& 4.50230000& 3.50770000& 2.15040000& 1.33999000 \\
   GPEAae& 6.0148080& 5.30672400& 4.49323200& 3.50516800& 2.15378100& 1.33998200 \\
   HPO& 6.00552336& 5.30591367& 4.49474956& 3.51336235& 2.15423400& 1.33652825 \\ 
      IAS& 5.9914000& 5.30850000& 4.51190000& 3.50210000& 2.16010000& 1.34000000 \\
      MVO& 6.0239402& 5.30601120& 4.49501130& 3.49602200& 2.15272610& 1.33995950 \\
      SOS& 6.0187800& 5.30344000& 4.49587000& 3.49896000& 2.15564000& 1.33996000 \\
     \hline
    \end{tabular*}
    }
    \caption{Comparison results for the cantilever design problem.}
    \label{Tab3}   
\end{table}

\subsubsection{Pressure vessel design problem} 

\begin{figure}[ht]
    \centering
    \includegraphics[width=0.8\textwidth]{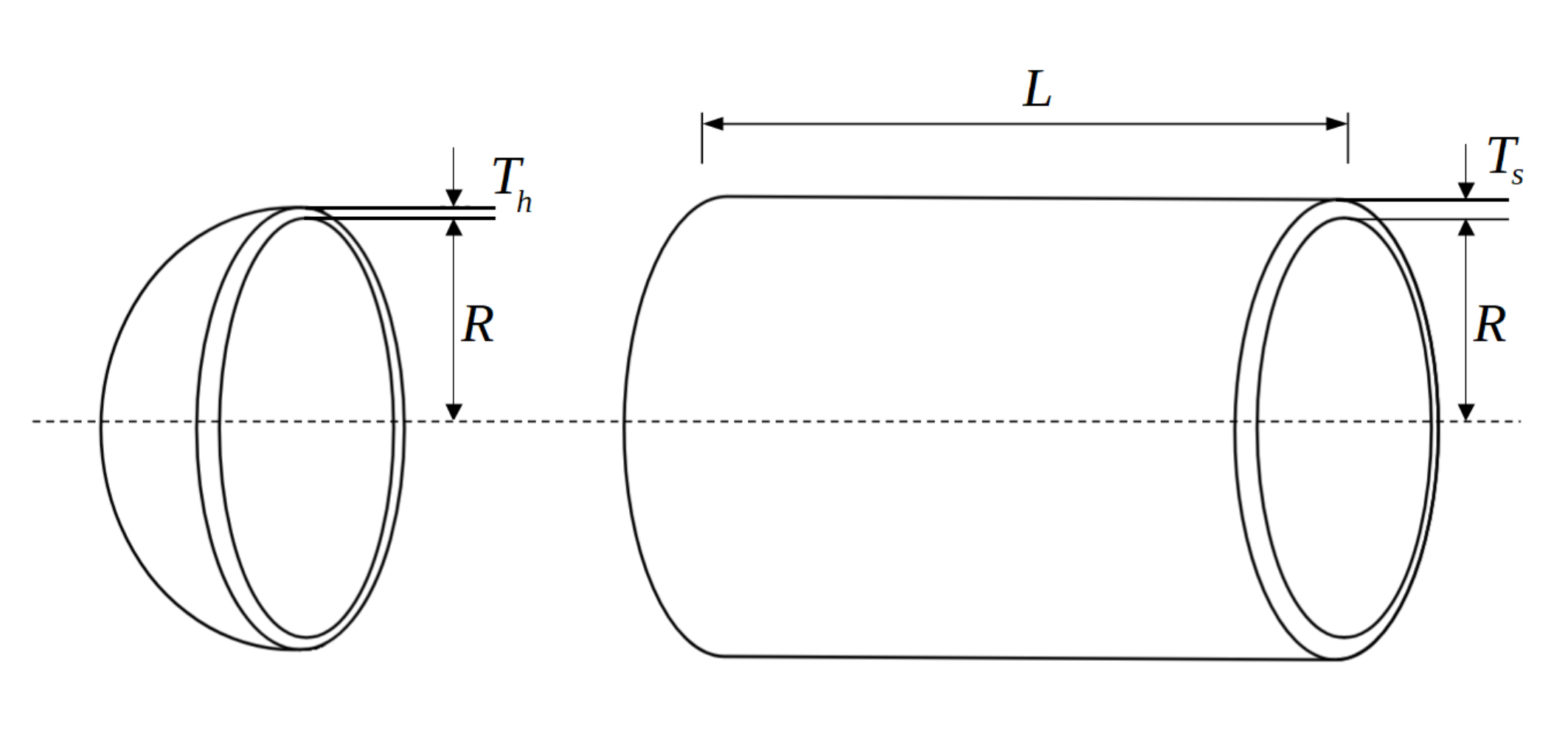}
    \caption{Pressure vessel design problem.}
    \label{figura3}
\end{figure}

The aim of this problem is to achieve cost minimization encompassing material, forming, and welding expenses associated with a cylindrical vessel, as illustrated in Figure \ref{figura3}. The vessel is enclosed at both ends with a hemispherical-shaped head. There are four key variables in this problem: $T_s$ ($x_1,$ thickness of the shell), $T_h$ ($x_2,$ thickness of the head), $R$ ($x_3,$ inner radius) and $L$ ($x_4,$ length of the cylindrical section of the vessel, not including the head).

This problem is governed by four constraints. The formulation of these constraints and the problem is articulated as follows \cite{Kan}:

\begin{align} \label{eq9}
    & \min \quad  f(x) = 0.6244x_1x_3x_4 + 1.7781x_2x_3^2 +3.1661x_1^2x_4 + 19.84x_1^2x_3  \\
    & \text{Subject to:}  \nonumber \\ 
    & g_1(x): -x_1 + 0.0193x_3\leq 0 \nonumber \\
    & g_2(x): -x_2 + 0.00954x_3\leq 0 \nonumber \\
    & g_3(x): -\pi x_3^2x_4 - \frac{4}{3}\pi x_3^3 +1296000 \leq 0 \nonumber \\ 
    & g_4(x): x_4 - 240\leq 0 \nonumber \\
    & 0 \leq x_1 \leq 99 \nonumber \\
    & 0 \leq x_2 \leq 99 \nonumber \\
    & 10 \leq x_3 \leq 200 \nonumber \\
    & 10 \leq x_4 \leq 200. \nonumber 
\end{align}

To solve problem (\ref{eq9}), we transform it into an unconstrained problem following strategy (\ref{eqOpt2}), with the number of constraints being $m=4.$ TAO algorithm yields the following solution: $x^+=(0.77873582, 0.38572842, 40.34900616, 199.59130975),$ with $f(x^+)=5888.6156573066.$
In accordance with the findings reported in the article \cite{Naru}, the comparative results with other algorithms are presented in Table \ref{Tab4}. Algorithms under consideration include AEO, BA, charged system search (CSS), genetic algorithm (GA), GPEAae, Gaussian quantum-behaved particle swarm optimization (G-QPSO), GWO, Harris hawks optimizer (HHO), HPO, hybrid particle swarm optimization (HPSO), mothflame optimization (MFO), sine cosine gray wolf optimizer
(SC-GWO), water evaporation optimization(WEO) and whale optimization algorithm (WOA). As inferred from the findings, the proposed TAO algorithm has demonstrated excellent performance in addressing this problem.

\begin{table}[ht]
{\small
    \noindent\begin{tabular*}{\textwidth}{l @{\extracolsep{\fill}}lllll}
    \toprule
    \multicolumn{1}{l}{\small Algorithm}& \multicolumn{4}{l}{\small Optimum variables}& \multicolumn{1}{l}{\small Optimum cost}\\
    \cline{2-5} 
        &$T_s$& $T_h$& $R$& $L$& \\
      TAO& 0.77873582& 0.38572842&  40.34900616& 199.59130975& 5888.6156573066 \\
      AEO& 0.8374205& 0.413937& 43.389597& 161.268592& 5994.50695 \\   
      BA& 0.812500& 0.437500& 42.098445& 176.636595& 6059.7143 \\
      CSS& 0.812500& 0.437500& 42.103624& 176.572656& 6059.0888 \\
      GA& 0.812500& 0.437500& 42.097398& 176.654050& 6059.9463 \\
     GPEAae& 0.812500& 0.437500& 42.098497& 176.635954& 6059.708025 \\ 
     G-QPSO& 0.812500& 0.437500& 42.0984& 176.6372& 6059.7208 \\
     GWO& 0.8125& 0.4345& 42.089181& 176.758731& 6051.5639 \\
     HHO& 0.81758383& 0.4072927& 42.09174576& 176.7196352& 6000.46259 \\
     HPO& 0.778168& 0.384649& 40.3196187& 200& 5885.33277 \\
     HPSO& 0.812500& 0.437500& 42.0984& 176.6366& 6059.7143 \\
     MFO& 0.8125& 0.4375& 42.098445& 176.636596& 6059.7143 \\
     SC-GWO& 0.8125& 0.4375& 42.0984& 176.63706& 6059.7179 \\
     WEO& 0.812500& 0.437500& 42.098444& 176.636622& 6059.71 \\
     WOA& 0.812500& 0.437500& 42.0982699& 176.638998& 6059.7410 \\
     \hline
    \end{tabular*}
    }
    \caption{Comparison of results for the pressure vessel design problem.}
    \label{Tab4}   
\end{table}

\subsubsection{Tension/compression spring design problem}
\begin{figure}[ht]
    \centering
    \includegraphics[width=0.6\textwidth]{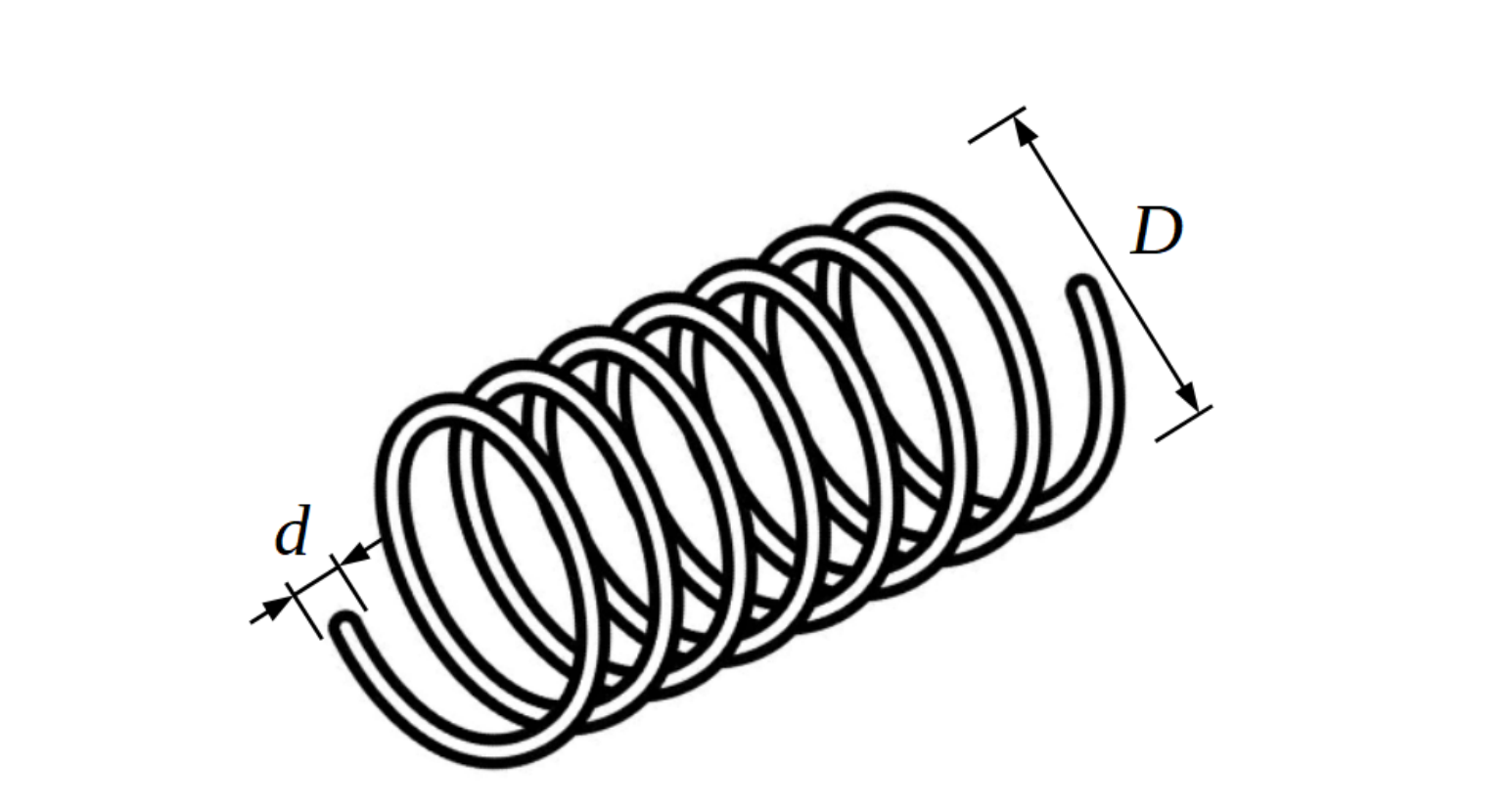}
    \caption{Tension/compression spring design.}
    \label{figura4}
\end{figure}

The engineering test problem employed in this study pertains to the optimization of tension/compression spring design. The objetive of this problem is to minimize the weight of a tension/compression spring associated with a spring characterized by three key parameters, namely the number of active loops ($N$), the average coil diameter ($D$), and the wire diameter ($d$). The geometric details of the spring and its associated parameters are visually presented in Figure \ref{figura4}.

The spring design problem is governed by a set of inequality constraints, formally expressed in equation (\ref{eq10}). 

\begin{align} \label{eq10}
    & \min \quad  f(x) = (x_3+2)x_2x_1  \\
    & \text{Subject to:}  \nonumber \\ 
    & g_1(x): 1 - \frac{x_2^3x_3}{71785x_1^4}\leq 0 \nonumber \\
    & g_2(x): \frac{4x_2^2-x_1x_2}{12566(x_2x_1^3-x_1^4)} + \frac{1}{5108x_1^2} - 1 \leq 0 \nonumber \\
    & g_3(x): 1 - \frac{140.45x_1}{x_2^2x_3} \leq 0 \nonumber \\ 
    & g_4(x): \frac{x_2+x_1}{1.5} - 1 \leq 0 \nonumber \\
    & 0.05 \leq x_1 \leq 2 \nonumber \\
    & 0.25 \leq x_2 \leq 1.3 \nonumber \\
    & 2.00 \leq x_2 \leq 15, \nonumber 
\end{align}
where $x_1 = d,$ $x_2 = D$ and $x_3=N.$

To solve problem (\ref{eq10}), we have applied the methodology described in approach (\ref{eqOpt2}), with the number of constraints in this problem being $m=4.$

\begin{table}[ht]
{\small
    \noindent\begin{tabular*}{\textwidth}{l @{\extracolsep{\fill}}llll}
    \toprule
    \multicolumn{1}{l}{\small Algorithm}& \multicolumn{3}{l}{\small Optimum variables}& \multicolumn{1}{l}{\small Optimum weigth}\\
    \cline{2-4} 
        &$N$& $D$& $d$&  \\
       TAO& 9.65579408& 0.38821894& 0.05296587& 0.0126943602 \\
       AEO& 10.879842& 0.361751& 0.051897& 0.0126662 \\   
       BA& 11.2885& 0.35673& 0.05169& 0.012665 \\
       COOT& 11.34038& 0.35584& 0.05165& 0.012665293 \\
       CPSO& 11.244543& 0.357644& 0.051728& 0.012674 \\
       GPEAae& 11.294000& 0.356631& 0.051685& 0.012665 \\
       GWO& 11.28885& 0.356737& 0.05169& 0.012666 \\
       HHO& 11.138859& 0.359305355& 0.051796393& 0.012665443 \\
       HPO& 11.21536452893& 0.3579796674& 0.0517414615& 0.012665282823723 \\
       SFS& 11.288966& 0.356717736& 0.051689061& 0.012665233 \\
	   SHO& 12.09550& 0.343751& 0.051144& 0.012674000 \\
	   SSA& 12.004032& 0.345215& 0.051207& 0.0126763 \\
       WCA& 11.30041& 0.356522& 0.05168& 0.012665 \\
       WEO& 11.294103& 0.356630& 0.051685& 0.012665 \\
       WOA& 12.004032& 0.345215& 0.051207& 0.0126763 \\
     \hline
    \end{tabular*}
    }
    \caption{Comparison of results for the Tension/compression spring problem.}
    \label{Tab5}   
\end{table}

Table \ref{Tab5} presents a comparative analysis of the results obtained using the TAO algorithm alongside other heuristic optimization methods in the context of this specific problem. This table provides comprehensive information concerning the values of each design variable and their respective optimization outcomes. 
Algorithms considered in the comparison of results are the following: AEO, BA, COOT, CPSO, GPEAae, GWO, HPO, HHO, stochastic fractal search (SFS), spotted hyena optimizer (SHO), salp swarm algorithm (SSA), water cycle algorithm (WCA), WEO and WOA. TAO algorithm produces the following solution: $x^+=(0.05296587 0.38821894 9.65579408),$ where $f(x^+)=0.0126943602.$

\subsubsection{Sustainable exploitation of a natural renewable resource}

In \cite{Her}, a marine resource exploited from 1999 to 2019 is studied. It is assumed that the resource's biomass, measured in tons, follows the Gordon-Schaefer equation, which is mathematically expressed as:

$$
B(t+1) = B(t) + rB(t)\left(1-\frac{B(t)}{K}\right) - h(t), t=1999, \dots,2019,
$$
where $B(t)$ represents the resource's biomass in year $t$, $r$ and $K$ are, respectively, the growth parameter and the carrying capacity of the fish population, and $h(t)$ is the catch made in year $t$. The catch is defined as $h(t)=qE(t)B(t),$ where $E(t)$ is the fishing effort, representing the intensity of human activities to extract fish, and $q$ is the catchability coefficient, representing the fraction of the population captured per unit of fishing effort. The unit of measurement for $E(t)$ is boats-year (and, consequently, the unit of measurement for $q$ is 1/boats-year).

Using the Bayesian estimation MCMC method, the distributions of $r,$ $K,$ and $q$ were calculated based on empirical data  \cite{Her}. The mean values of the estimated parameters are: $r=0.7534,$ $K=2.7399\times10^4$ Tn, and $q=0.01081$ 1/boats-year. In 2019 (the initial time for our analysis), the biomass of the resource was $1.0836\times10^4$ Tn. Thus, we have the following model, with the estimated parameters, for the dynamics of the biomass of the resource subject to harvesting:

\begin{equation}\label{eq11}
\begin{cases}
B(t+1) = B(t) + rB(t)\left(1-\frac{B(t)}{K}\right) - qE(t)B(t) \\
B(0) = 1.8836\times10^4 \text{ Tn}.
\end{cases}
\end{equation}

Now, we consider the problem of finding a harvesting strategy that allows for the exploitation of the fishery resource while maximizing the gains for the entity exploiting that resource. In other words, we seek strategies that maximize the catch benefits while simultaneously protecting the species over a time period $[2019, 2019+T],$ where $T$ is known as the time horizon. To model the catch benefits, we can observe that if $p$ is the price per unit of fishing effort, then the profit at time $t$ is $pqE(t)B(t).$ On the other hand, if $c$ is the unit cost per unit of fishing effort, then the losses at time $t$ are $cE(t),$ and the rent in year $t$ is given by $pqE(t)B(t)-cE(t).$ Then, the benefit obtained $T$ years after 2019 is modeled as the following sum:

\begin{equation}\label{eq12}
	\displaystyle\sum_{t=0}^{T-1}\rho^t\left(pqE(t)B(t)-cE(t)\right),
\end{equation}
where $\rho\in(0,1].$

In order to protect the exploited species, strategies are sought to maximize its biomass at the final time $T,$ i.e., the quantity $B(T).$ Therefore, the objective is to maximize the following functional:

\begin{equation}\label{eq13}
    J(E) = \displaystyle\sum_{t=0}^{T-1}\rho^t\left(pqE(t)B(t)-cE(t)\right)+B(T),
\end{equation}

where $\rho$ is a discount factor chosen in $[0,1).$

Interpretation of (\ref{eq13}) is as follows: as time passes, less weight (or importance) is given to the benefit obtained by the entity that catches fish, and more weight is given to the protection of the resource, because, as $t$ increases, $\rho$ decreases. Furthermore, for biological reasons, it may be required that in each year, $t,$ $B(t) \in [B_u, B_s],$ where $B_u$ and $B_s$ are minimum and maximum biomass threshold. 

Thus, the following optimization problem is formulated:

\begin{align} \label{eq14}
    & \max \quad  J(E) = \displaystyle\sum_{t=0}^{T-1}\rho^t\left(pqE(t)B(t)-cE(t)\right)+B(T)  \\
    & \text{Subject to:}  \nonumber \\ 
    & \begin{cases}
        B(t+1) = B(t) + rB(t)\left(1-\frac{B(t)}{K}\right) - qE(t)B(t) \\
        B(0) = 1.8836\times10^4 \text{ Tn}
\end{cases} \nonumber \\
    & g(E): B_u - B(t) \leq 0 \nonumber \\
    & E_{min} \leq E(t) \leq E_{max}, t=2019, \dots, 2019+T, \nonumber 
\end{align}
where $E_{min}$ and $E_{max}$ are the minimum and maximum values of fishing effort.

\begin{table}[ht]
{\small
    \noindent\begin{tabular*}{\textwidth}{l @{\extracolsep{\fill}}ll}
    \toprule
       Parameter& Description& Value \\
       \hline
       $B(2019)$& Initial biomass& $1.0836\times10^4$ Tn \\
       $r$& Growth parameter& 0.7534 \\
       $K$& Carrying capacity& $2.7399\times10^4$ Tn \\
       $q$& Catchability coefficient& $0.01081$/boats-year \\
       $E_{min}$& Minimum fishing effort& 0 boats-year \\ 
       $E_{max}$& Maximum fishing effort& 41 boats-year \\
       $p$& Price per unit of fishing effort& \$ 6 000 \\
       $c$& Unit cost per unit of fishing effort& \$ 3 070 \\
       $T$& Time horizon& 30 \\
       $\rho$& Discount factor& $[0,1)$ \\ 
       $B_u$& Minimum biomass threshold& $1.0836\times10^4$ Tn \\
       $B_s$& Maximum biomass threshold& $2.2\times10^4$ Tn \\          
     \hline
    \end{tabular*}
    }
    \caption{Parameters used in problem (\ref{eq14}).}
    \label{Tab6}   
\end{table}

We used the parameters in Table \ref{Tab6} to conduct a 30-year sustainability study. For simplicity, instead of considering the interval $[2019, 2019+T]$, we have used the interval $[1, T]$. We assume a scenario in which it is required that the biomass of the resource does not decrease below the initial biomass, that is, that of 2019. In addition, weight is assigned to the profits obtained, which leads us to consider a high value for $ \rho$, specifically, $\rho=0.9$. To solve the optimization problem (\ref{eq14}), we follow the strategy (\ref{eqOpt1}).

\begin{figure}[ht]
    \centering
    \includegraphics[width=0.9\textwidth]{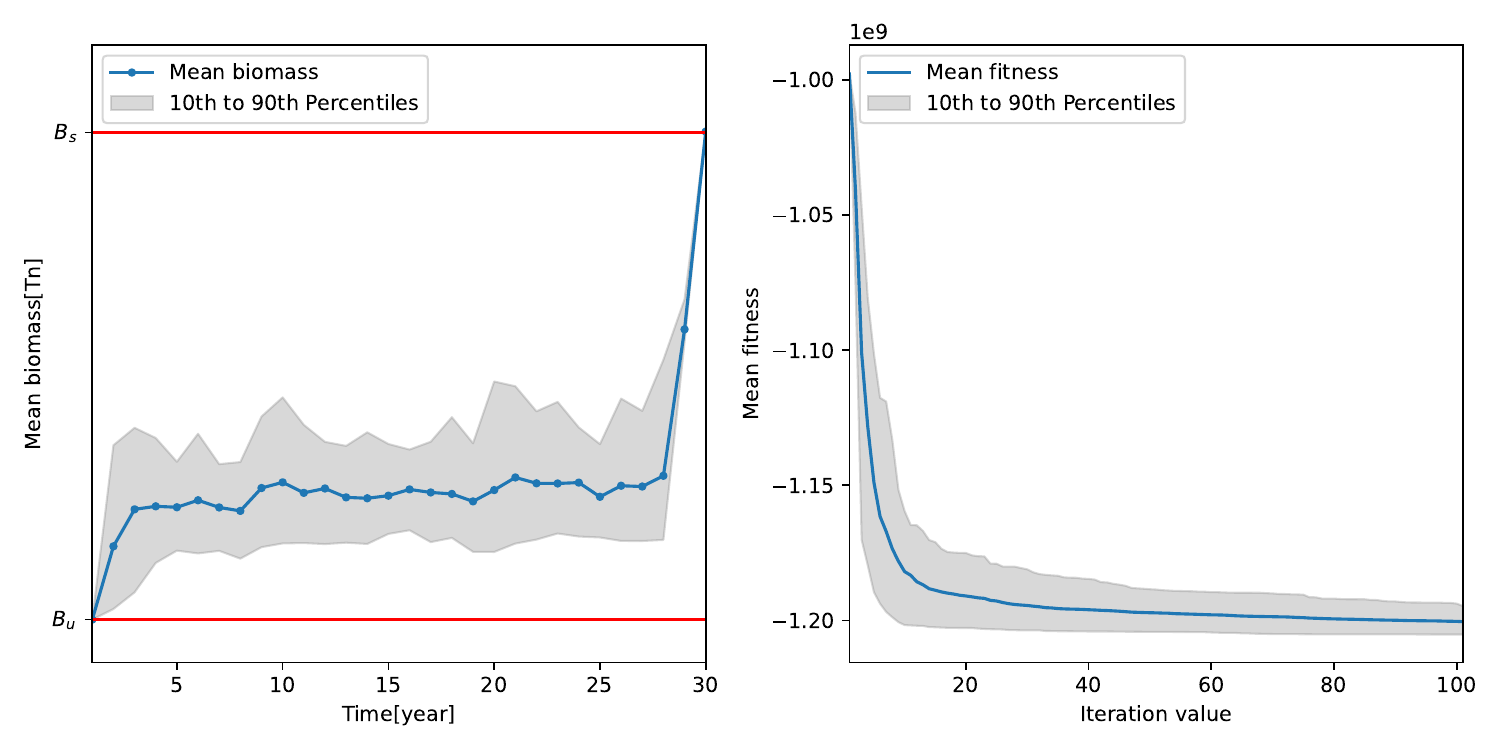}
    \caption{Sustainable exploitation renewable resource: mean biomass (left) and convergence plot (right).}
    \label{figura5}
\end{figure}

\begin{figure}[ht]
    \centering
    \includegraphics[width=0.8\textwidth]{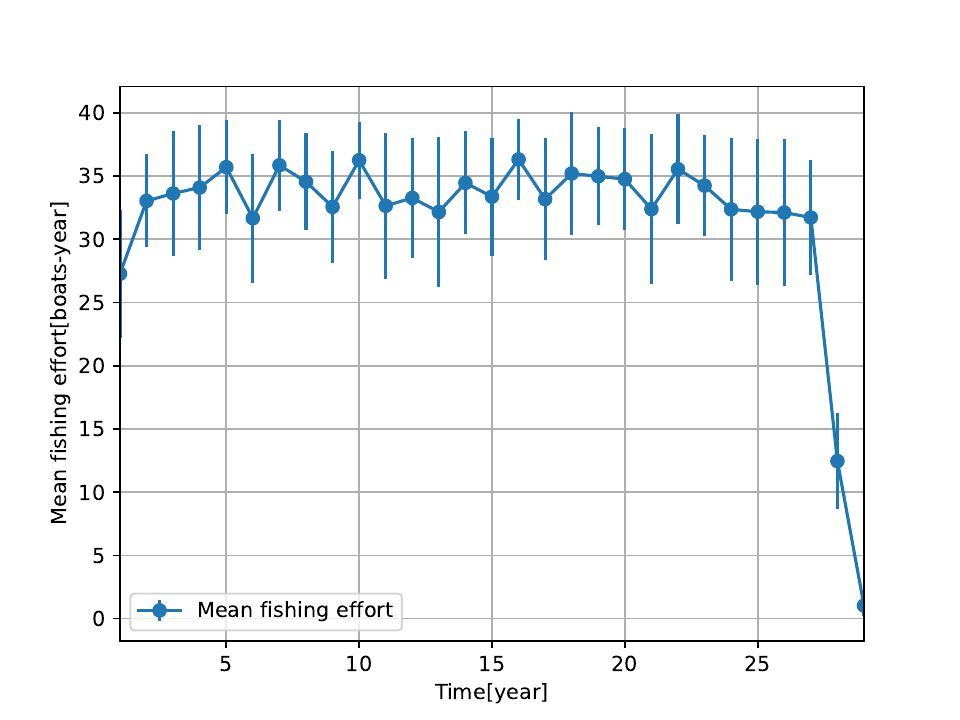}
    \caption{Mean fishing effort[boats-year] dynamic.}
    \label{figura6}
\end{figure}

We realized 100 simulations, with population size of 50 cells and 100 iterations. The average biomass dynamics (see the right side of Figure \ref{figura5}) show that it is possible to exploit the resource sustainably, in the previously defined sense, with an average annual fishing effort such as that shown in Figure \ref{figura6}.

Various sustainability scenarios can be analyzed by considering different values of the parameter $\rho$ and different time horizons, $T.$ For example, varying the value of $T$ allows us to assess the short, medium, and long-term impact of resource exploitation. 
Moreover, if $T$ is fixed, we can examine the impact of exploitation based on the emphasis placed on profits versus the preservation of the exploited resource.

\section{Conclusiones}

In this article, we present a new metaheuristic algorithm called Tumor Angiogenesis Optimizer (TAO). This algorithm is inspired by the morphogenetic cellular movements exhibited by ECs during the complex process of tumor angiogenesis. According to the results obtained from test functions, TAO outperforms the PSO standard. Furthermore, our algorithm has demonstrated high performance in real-world constrained optimization problems. In particular, in the context of cantilever beam design problems and pressure vessel design problems, it outperformed several well-established metaheuristics, while in the tension/compression spring design problem it had acceptable performance. In this article, we have also successfully addressed a sustainability problem related to the exploitation of a renewable resource.

For future work, the proposed algorithm can be used in different fields of study. In addition, we will delve into the mathematical modeling of the cellular angiogenesis process, focusing specifically on the phenomena of cell emergence and death, which could make the TAO algorithm more efficient. Furthermore, we aim to develop a version of our algorithm adapted to address multi-objective optimization problems.

\end{document}